\documentclass[12pt,notitlepage]{article}
\usepackage{amsmath}
\usepackage{amssymb}
\usepackage{amscd}
\usepackage{amsfonts}
\usepackage{amsthm}
 \theoremstyle{definition}
 \theoremstyle{plain}
  \newtheorem{theorem}{Theorem}
  \newtheorem{prop}[theorem]{Proposition}
  \newtheorem{lemma}[theorem]{Lemma}
  \newtheorem{cor}[theorem]{Corollary}
 \theoremstyle{remark}
  \newtheorem*{rem}{\bf Remark}

\newcommand{\Dj}{\hbox to 8pt{\raisebox{.4\height}{-}\hss D}}
\newcommand{\dJ}{\,{\hbox to 4pt{\raisebox{.75\height}{-}\hss d}}}

\newcommand{\ac}{\ensuremath{\mathcal A}}

\newcommand{\hc}{\ensuremath{\mathcal H}}

\newcommand{\eqdef}{\stackrel{\rm def}{=}}
\begin{document}

\title{\begin{flushright}
\small\bf ITEP-TH-65/01\\
\end{flushright}
\vspace{10mm}
\bf\large Hopf-type Cyclic Cohomology\\
via Karoubi Operator}
\author{\em G.~I.~Sharygin${}^{1,2}$
\thanks{The author was partially supported by the
Russian Foundation for Basic Research (grant
no.~01-01-00546).}}
\date{September-December 2001}
\maketitle
\begin{center}
${}^{1}${\em ITEP, Russia, 117218, Moscow, B. Cheremushkinskaya, 25}\\
${}^2${\em Kolmogorov College, Russia, 121357, Moscow, Kremenchugskaya, 11}
\end{center}

\begin{abstract}
In this paper we propose still another approach to the Hopf-type
cohomologies of a Hopf algebrae \hc, based on the notion of the universal
differential calculus on \hc. Few remarks, concerning the possible generalizations
and applications of this approach are made.
\end{abstract}
\section{Introduction}
In the original papers of A.~Connes and H.~Moscovici (see
\cite{{CM1}, {CM2}, {CM3}}) the explicit structure of cyclic
module defining the so-called Hopf-type cyclic cohomologies of a
Hopf algebra was given. Later, in his paper M.~Crainic showed that
this cyclic module could be obtained as the space of coinvariants
of the Hopf algebra's action on some other cyclic module. Some
further generalizations and developments, we know of, were made in
the papers \cite{theil}, \cite{iran1}, \cite{iran2}.

The purpose of this paper is to describe the Hopf-type cohomologies
of a Hopf algebra in the terms of a subcomodule of the so-called
algebra of non-commutative differential forms, associated with the
Hopf algebra. It turns out, that for any "modular pair in
involution", $(\delta,\sigma)$ one can associate a subcomplex of
this differential algebra, stable under Karoubi operator $\kappa$
(see papers \cite{{CQ1},{CQ2}} and \cite{Kar}) or its twisted
version $\kappa_\xi$, see \eqref{eq3.9} (also \cite{tuset}).

In addition to giving a new point of view on this homology theory,
this approach seems to have some virtues of its own. For instance,
one can try to define some similar sort of cyclic cohomologies, when
modular pair is substituted for a more general object. Besides
this, it can be used to establish bridges between this cyclic
cohomology theory and the Hopf-Galois theory, developed in the
papers of T.~Brzezinski, M.~\Dj ur\dJ evi\v{c}, P.~Hajac, S.~Majid
\cite{{Take},{Durd1},{Durd2},{Durd3},{Maj-Brz},{Haj}} and others. Only few
remarks, concerning this subject are made here, since we postpone
deeper discussions to a paper to follow.

Let's, first of all recall the construction of the Hopf-type
cohomologies, due to A.~Connes and H.~Moscovici. Here and below \hc\
will denote a Hopf algebra over a field of
characteristic 0 ($\mathbb C$ is our main example). Let $m,\ \Delta,\ 1,\ \epsilon$ and $S$ be the
multiplication, comultiplication (or diagonal), unit, counit and
antipode of \hc\ respectively. Below we shall usually miss $m$ in
our formulae, and use the standard (Sweedler, \cite{Sweed}) notation to write down the
diagonal: $$ \Delta(h)=\sum h_{(1)}\otimes h_{(2)}. $$

One says, that $(\sigma,\,\delta)$, where $\sigma$ is a group-like
element in \hc\ and $\delta:\hc\to\mathbb{C}$ an algebraic
character, and $\delta(\sigma)=1$, is a modular pair in
involution, if
\begin{align}
\label{eq0}
S_\delta^2(h)=&\sigma\/h\/\sigma^{-1},\ h\in\hc,\\
\intertext{where} S_\delta(h)=&\sum\delta(h_{(1)})S(h_{(2)}).
\end{align}
This is equivalent to $(\sigma S_\delta)^2=1$.

Given a modular pair in involution one can define the (co)cyclic module
$\hc^\sharp_{(\delta,\sigma)}$. Recall, that cyclic category is self-dual,
hence it is not necessary to distinguish very carefully between cyclic and cocyclic objects. So, one
puts: $(\hc^\sharp_{(\delta,\sigma)})_n=\hc^{\otimes n}$, and the cyclic
structure maps are defined as follows:
\begin{align}\delta_i:\hc^\sharp_n\to\hc^\sharp_{n+1},\quad&i=0,\dots,n+1\\
\sigma_i:\hc^\sharp_n\to\hc^\sharp_{n-1},\quad&i=1,\dots,n,\\
\tau_n:\hc^\sharp_n\to\hc^\sharp_n.&
\end{align}
are given by
\begin{gather}
\delta_i(h_1,\dots,h_n)=
\begin{cases}
(1,h_1,\dots,h_n), &i=0,\\ (h_1,\dots,\Delta(h_i),\dots,h_n),
&1\leq i\leq n,\\ (h_1,\dots,h_n,\sigma), &i=n+1;
\end{cases}\\
\sigma_i(h_1,\dots, h_n)=\epsilon(h_i)(h_1,\dots,h_{i-1},\hat
h_i,h_{i+1},\dots,h_n),\ 1\leq i\leq n,\\
\tau_n(h_1,h_2,\dots,h_n)=S_\delta(h_1)\cdot(h_2,\dots,h_n,\sigma).
\end{gather}
Here in the last formula we assume, that \hc\ acts on its own
tensor power as follows: $$
h\cdot(h_1,\dots,h_n)=(h_{(1)}h_1,\dots,h_{(n)}h_n).$$

For any cocyclic module one can define its cyclic, negative cyclic
and periodic cyclic cohomology. To this end one has to consider the
cyclic, negative and periodic complexes respectively (see, for
example the book of Loday \cite{Loday} and paper \cite{CQ2}). For instance,
periodic cohomology are defined by the following super-complex:
\begin{equation}
CP_i=\prod_{n\equiv i({\rm mod\/}2)}\hc^\sharp_n,\ i=0,1;
\end{equation}
equipped with differentials $b:\hc^\sharp_n\to\hc^\sharp_{n+1}$
and $B:\hc^\sharp_n\to\hc^\sharp_{n-1}$, defined as follows
\begin{align}
\label{eq0'} b=&\sum_i(-1)^i\delta_i,\\ \label{eq0"}
B=&N\circ(1-\tau_{n+1})\circ\tilde\sigma_0,\\
\intertext{where}\label{eq0'"}\tilde\sigma_0=&\sigma_n\circ\tau_n.
\end{align}
Recall, that we deal with cocyclic module here, so the usual
formulae for differentials in mixed complexes, associated with
such an object, are inverted.

\section{Special case: $(\delta,\sigma)=(\epsilon,1)$}
Let $\mathcal H$ be the given (unital) Hopf algebra. We shall
denote by $\Omega(\mathcal H)$ the universal unital differential
graded algebra, generated by $\mathcal H$. Recall, that
\begin{gather}
\label{eq1} \Omega(\mathcal H) = \bigoplus_{n\ge 0}
\Omega_n(\mathcal H);\\ \label{eq2}\Omega_0(\mathcal H)=\mathcal
H;\quad\Omega_1(\mathcal H)=\ker(m:\mathcal H\otimes\mathcal
H\to\mathcal H);\\ \label{eq2} \Omega_n(\mathcal
H)=\underbrace{\Omega_1(\mathcal H)\otimes_{\mathcal H}
\Omega_1(\mathcal H)\otimes_{\mathcal H}\dots\otimes_{\mathcal H}
\Omega_1(\mathcal H)}_n.
\end{gather}
The differential $d:\mathcal H\to\Omega_1(\mathcal H)$ is given by
$$ d(x)=1\otimes x-x\otimes1,$$ and one can prove that any element
$\theta$ in $\Omega_n(\mathcal H)$ can be written down in the form
$$\theta=\sum_i a_0^i da_1^i da_2^i\dots da_n^i,\quad
a_j^i\in\hc$$ Now it's clear, that \begin{equation} \label{eq4}d\theta=\sum_i da_0^i da_1^i
da_2^i\dots da_n^i.\end{equation}

So far, the coalgebra structure hasn't yet come to the scene.
In effect, one can define the universal differential algebra,
associated to any unital algebra \ac\ in precisely same way. But
now, since \hc\ is in fact a Hopf algebra, one can define left-
and right-coactions of \hc\ on $\Omega(\hc)$. Namely, put
\begin{align}
\label{eq5}
\Delta_R(\theta)&=\sum_i a_{0,(1)}^i da_{1,(1)}^i
da_{2,(1)}^i\dots da_{n,(1)}^i\otimes a_{0,(2)}^i a_{1,(2)}^i
a_{2,(2)}^i\dots a_{n,(2)}^i,\\
\intertext { and}
\label{eq6}
\Delta_L(\theta)&=\sum_i a_{0,(1)}^i a_{1,(1)}^i a_{2,(1)}^i\dots
a_{n,(1)}^i\otimes a_{0,(2)}^i da_{1,(2)}^i da_{2,(2)}^i\dots
da_{n,(2)}^i.
\end{align}
 The fact, that these formulae really determine a well-defined
maps follows from the universal properties of $\Omega(\hc)$.

Moreover, \eqref{eq5} and \eqref{eq6} define the right- and left-
Hopf-comodule algebra structures on $\Omega(\hc)$.  That is, the map
$\Delta_R:\Omega(\hc)\to\Omega(\hc)\otimes\hc$ is an algebra morphism, and similarly
for $\Delta_L$. In particular, $\Omega(\hc)$ is a left and right Hopf module
over the Hopf algebra $\hc=\Omega_0(\hc)$. In a general case
conditions that formulas \eqref{eq5} and \eqref{eq6} define such structures impose an
additional restrictions on the structure of a given differential
calculus $\Omega'(\hc)$. Differential calculi verifying these restrictions are called {\em
bicovariant\/}. This matter is accurately explained in \cite{Wor}, where
the general definition of a bicovariant differential calculus on a Hopf
algebra is given.

Observe, see \eqref{eq4}, that both maps \eqref{eq5} and
\eqref{eq6} preserve differential.  Hence, in particular, the
subspaces of left and right coinvariants are differential graded
subalgebrae in $\Omega(\hc)$. Let's describe explicitly the
structure of these subalgebrae. For instance, take $\Omega^R(\hc)\eqdef\Omega(\hc)^{{\rm
co}\hc}$.

It is shown in \cite{Wor} and \cite{{Durd1},{Durd2}} that the map
\begin{equation}
\label{eq7'} \pi^R:\hc\to\Omega_1^R(\hc),\quad h\mapsto
da_{(1)}\cdot S(a_{(2)})
\end{equation}
identifies the space $\Omega_1^R(\hc)$ with $\ker\epsilon$.
Moreover, one can show, that
\begin{equation}\label{eq6'}\Omega(\hc)\cong\Omega^R(\hc)\otimes\hc,\quad\mbox{and
}\ \Omega^R_n(\hc)\cong\Omega^R_1(\hc)^{\otimes
n}\cong(\ker\epsilon)^{\otimes n}.\end{equation} In terms of these
isomorphisms, one can write down the differential and left and
right actions of \hc\ on the bimodule $\Omega(\hc)$ as follows:
\begin{gather}
dh=dh_{(1)}\cdot S(h_{(2)})h_{(3)}=\pi^R(h_{(1)})\otimes h_{(2)},\\
d\pi^R(h)=dh_{(1)}dS(h_{(2)})=-dh_{(1)}S(h_{(2)})dh_{(3)}S(h_{(4)})=\pi^R(h_{(1)})\otimes\pi^R(h_{(2)}),\\
\label{eq7}a\cdot\pi^R(h)=a_{(1)}\pi^R(h)S(a_{(2)})a_{(3)}=\pi^R(a_{(1)}h-\epsilon(h)a_{(1)})\otimes a_{(2)},
\end{gather}
and the right action of \hc\ on
$\Omega_1(\hc)=\Omega^R(\hc)\otimes\hc$ is trivial. In particular,
if $h$ is in $\ker\epsilon$, then formula \eqref{eq7} gives the
following description of the left action of \hc\ on bimodule
$\Omega_1(\hc)\cong\ker\epsilon\otimes\hc$:
\begin{equation}
\label{eq8}
 a\cdot (h_1\otimes h_0)=a_{(1)}h_1\otimes
a_{(2)}h_0,\quad a,h_0\in\hc,\ h_1\in\ker\epsilon.
\end{equation}
This formula extends in a natural way to the $n-$th graded
component of $\Omega(\hc),\
\Omega_n(\hc)\cong(\ker\epsilon)^{\otimes n}\otimes\hc$, namely
\begin{equation}
\label{eq9}
 a\cdot(h_1\otimes h_2\otimes\dots\otimes h_n\otimes
h)= a_{(1)}h_1\otimes a_{(2)}h_2\otimes\dots\otimes
a_{(n)}h_n\otimes a_{(n+1)}h,
\end{equation}
where $h_1,\dots h_n\in\ker\epsilon,\ h,a\in\hc.$

Let now $b$, and $\kappa$ be the usual Hochschild differential and
Karoubi operator on $\Omega(\hc)$, defined as follows (see
\cite{CQ2}):
\begin{align}
\label{eqnew1}
b(\omega da)&= (-1)^{|\omega|}(\omega a - a\omega),\\
\intertext{where $a\in\hc,\ \omega\in\Omega(\hc)$, and} \label{eqnew2}\kappa&=1-
bd -db.
\end{align}
Explicitly one can show, that $$\kappa(\omega\,
da)=(-1)^{|\omega|}da\,\omega.$$

This is the usual way these operators are introduced. For our
purposes it woud be useful to consider a little bit different
operators, $b',\ \kappa'$:
\begin{align}
\label{eq10}
b'(da\omega)&= a\omega - \omega a,\\
\label{eq11'}
\kappa'(da\omega)&=b'd+db'-1.\\
\intertext{Explicitly}
\label{eq10'}
\kappa'(da\omega)&=(-1)^{|\omega|}\omega da.
\end{align}

Let $B'=\sum{\kappa'}^i\circ d$. Then $B'$ corresponds to the
operator $B$ from \cite{CQ2}. and these operators verify all the
usual properties of $b$ and $\kappa$ and $B$, see \S3 of
\cite{CQ2}. This can be proven by a slight modification of the
reasoning used in the quoted paper in the usual setting. Hence, we
conclude,that $b'$ and $B'$ induce the structure of cyclic module
on $\Omega(\hc)$. Moreover, formulae \eqref{eq10} and
\eqref{eq10'} show that
\begin{align}
\label{eq11"}
\kappa'&=\kappa^{-1},\\
\intertext{and}
\label{eq11"'}
b'&=\pm b\circ\kappa'
\end{align}
The following theorem is the main result of this section.

\begin{theorem}
\label{theo1}
Let \hc\ be a Hopf algebra and the modular pair $(1,\epsilon)$ is
such, that the conditions of \cite{CM2} are satisfied (i.e.
$S^2=1$) then\ 
$\Omega^R(\hc)$ is a differential graded subalgebra in
$\Omega(\hc)$, stable under $b$ and Karoubi operator $\kappa$\
(and, $b'$ and $\kappa'$) and hence is a mixed subcomplex in
$(\Omega(\hc), b, B)$. The same is true about $\Omega^L(\hc)$.\
Moreover, periodic cohomologies of
$\Omega^R(\hc)$\ 
with the mixed complex structure, induced from $(\Omega(\hc), b',
B')$, are naturally isomorphic to the periodic Hopf-type
cohomologies $HP_{\epsilon,1}^{*}(\hc)$ of the Hopf algebra
\hc.
\end{theorem}
\begin{proof}
The fact, that $\Omega^R(\hc)$ (and $\Omega^L(\hc)$ as well) is a
differential graded subalgebra in $\Omega(\hc)$ follows directly
from the discussion above. Now we shall prove the second statement
of this theorem. We shall confine our attention to
$\Omega^R(\hc)$. (In the case of $\Omega^L(\hc)$ reasoning is
absolutely similar.)

First of all let's note, that it's enough to prove the stability
of $\Omega^R(\hc)$ under the action of $b'$ and $\kappa'$ (just
look at formulae \eqref{eq11"} and \eqref{eq11"'}). So, let's start
with  proving, that $\Omega^R(\hc)$ is stable under $b'$. To this
end we shall directly compute the image of an element
$\omega\in\Omega^R(\hc)$ under $b'$. First, let $\omega$ belong to
$\Omega_1^R(\hc)$. We compute:
\begin{equation}
\begin{split}
b'(\omega)&=b'(da_{(1)}S(a_{(2)}))
=a_{(1)}S(a_{(2)})-S(a_{(2)})a_{(1)}\\&=\epsilon(a)\cdot1-S(a_{(2)})S^2(a_{(1)})
=\epsilon(a)\cdot1-
S\Bigl(S(a_{(1)})a_{(2)}\Bigr)\\&=\epsilon(a)\cdot1-\epsilon(a)\cdot1
=0,
\end{split}
\end{equation}
which is, of course, a right-coinvariant element. Here we've used
formula \eqref{eq10} and the possibility to represent any element
in $\Omega_1^R(\hc)$ as the image of some $a\in\hc$ under the map
$\pi^R$ from \eqref{eq7'}.

Now, if the element $\omega$ belongs to $\Omega^R_n(\hc)$ we use
the identification \eqref{eq6'} and formula \eqref{eq9} to compute
\begin{equation}
\label{eq12}
\begin{split}
b'(\omega)&=b'\bigl(\pi^R(a_1)\otimes\pi^R(a_2)\otimes\dots\otimes\pi^R(a_n)\bigr)\\
&=b'(d(a_{1,(1)})S(a_{1,(2)})\omega')=a_{1,(1)}S(a_{1,(2)})\omega'-S(a_{1,(2)})\omega'a_{1,(1)}\\
&=-\Bigl(\pi^R(S(a_{1,(n+1)})a_2)\otimes\dots\otimes\pi^R(S(a_{1,(3)})a_n)\Bigr)\,(S(a_{1,(2)})a_{1,(1)})\\
&=-\Bigl(\pi^R(S(a_{1,(n+1)})a_2)\otimes\dots\otimes\pi^R(S(a_{1,(3)})a_n)\Bigr)\,(S(a_{1,(2)})S^2(a_{1,(1)}))\\
&=-\Bigl(\pi^R(S(a_{1,(n+1)})a_2)\otimes\dots\otimes\pi^R(S(a_{1,(3)})a_n)\Bigr)\,\Bigl(S(S(a_{1,(1)})a_{1,(2)})\Bigr)\\
&=-\pi^R(S(a_{1,(n-1)})a_2)\otimes\pi^R(S(a_{1,(n)})a_3)\otimes\dots\otimes\pi^R(S(a_{1,(1)})a_n).
\end{split}
\end{equation}
Here $a_i\in\ker\epsilon,\ i=1,\dots,n$, and we denote for brevity
$\omega'=\pi^R(a_2)\otimes\pi^R(a_3)\otimes\dots\otimes\pi^R(a_n)$.
Clearly, $b'(\omega)$ lies in $\Omega_{n-1}^R(\hc)$.

Of course, since the inverse of Karoubi operator $\kappa'$ is
written down in terms of $d$ and $b'$, (see \eqref{eq11'}) one can
conclude, that $\Omega^R(\hc)$ is stable under its action. So, the
second statement of our theorem is proved.

However below we shall need the explicit
formula for this operator written down in terms of identification \eqref{eq6'}.
We use \eqref{eq10'} and find (we stick to the notation explained
after \eqref{eq11'}):
\begin{multline}
\label{eq13}
\begin{split}
&\kappa'(da_{1,(1)}S(a_{1,(2)})\omega')=(-1)^{|\omega'|}S(a_{1,(2)}\omega'\,da_{1,(1)}\\
&=(-1)^{|\omega'|}\Bigl(\pi^R(S(a_{1,(n+1)})a_2)\otimes\dots\otimes\pi^R(S(a_{1,(3)})a_n)\Bigr)\,(S(a_{1,(2)})da_{1,(1)})\\
&=(-1)^{|\omega'|}\Bigl(\pi^R(S(a_{1,(n+1)})a_2)\otimes\dots\otimes\pi^R(S(a_{1,(3)})a_n)\Bigr)\,(S(a_{1,(2)})dS^2(a_{1,(1)}))\\
&=(-1)^{|\omega'|+1}\Bigl(\pi^R(S(a_{1,(n)})a_2)\otimes\dots\otimes\pi^R(S(a_{1,(2)})a_n)\Bigr)\,\pi^R(S(a_{1,(1)}))\\
&=(-1)^n\pi^R(S(a_{1,(n)})a_2)\otimes\dots\otimes\pi^R(S(a_{1,(2)})a_n)\otimes\pi^R(S(a_{1,(1)})-\epsilon(a_{,(1)})).
\end{split}
\end{multline}
In other words, this can be written down as
\begin{equation}
\begin{split}
\kappa':\ker\epsilon^{\otimes n}&\to\ker\epsilon^{\otimes n},\\
(h_1,h_2,\dots,h_n)&\mapsto (-1)^n{\rm
proj}'\bigl(S(a_1)\cdot(h_2,\dots,h_n,1)\bigr).
\end{split}
\end{equation}
Here ${\rm proj}'$ denotes the standard projection ${\rm proj}':
\hc^{\otimes n}\to\ker\epsilon^{\otimes n}$, which sends each
component $h_i$ to $h_i-\epsilon(h_i)$.

Now, we believe, the similarity between these formulae and the
structure of cyclic module, introduced by A.~Connes and
H.~Moscovici is conspicuous. For instance, the cyclic operator
$\tau$ of this module is given by $$
\tau(h_1,h_2,\dots,h_n)=S(h_1)(h_2,\dots,h_n,1),$$ that is it
coincides with $\kappa'$, up to the sign nd projection on the
kernel of counit.

In the view of this observation, let's finally show, that the
cohomology of the induced sub mixed complex $(\Omega^R(\hc), b',
B')$ coincide with the Hopf-type cohomologies of Connes and
Moscovici.

To this end we first consider the cyclic object
$\hc^\sharp_{(\epsilon,1)}$, defined in \cite{CM2} and
\cite{Crai}, see section 1 above. Let $\tilde b, \tilde B$, be the
differentials \eqref{eq0'}, \eqref{eq0'"} and \eqref{eq0'"} (we
use tilde here to distinguish these maps from the Cuntz-Quillen's
operators on $\Omega(\hc)$, introduced above).

Recall, that (in this special case)
\begin{align}
\delta&=\sum (-1)^i\delta_i,\\ \intertext{where}
\delta_i(h_1,h_2,\dots,h_n)&=\begin{cases}
                               (1,h_1,h_2,\dots,h_n), &i=0,\\
                               (h_1,\dots,\Delta(h_i),\dots,h_n),
                               &1\leq i\leq n,\\
                               (h_1,h_2,\dots,h_n,1), &i=n+1,
                            \end{cases}\\
\tilde B&=N\circ\tilde\sigma_0\circ(1-\tau_{n+1})\\
\intertext{where} \notag
\tilde\sigma_0(h_1,h_2,\dots,h_n)&=S(h_1)\cdot(h_2,\dots,h_n),\\
\notag N&=\sum_{i=0}^n (-1)^{in}\tau_n^i.
\end{align}

Now, consider a slightly different mixed complex
$(\tilde\hc^\sharp_{(\epsilon,1)}, \tilde b', \tilde B')$, where
$(\tilde\hc^\sharp_{(\epsilon,1)})_n = \ker\epsilon^{\otimes n}$ and $\tilde
b'={\rm proj}'\circ \tilde b,\ \tilde B'={\rm proj}'\circ\tilde
B$.
\begin{lemma}
The natural projection from $\hc^\sharp_{(\epsilon,1)}$ to
$\tilde\hc^\sharp_{(\epsilon,1)}$ induces isomorphism on cyclic
cohomology.
\end{lemma}
\begin{proof}
This is a direct consequence of the fact that this projection
yields an isomorphism of the Hochschild homologies of these two
compexes (i.e. their homologies with respect to the differentials
$\tilde b$ and $\tilde b'$), which is a standard fact of homology
algebra (in fact, the latter complex is just the normalization of
the former one with respect to the degeneracy operators
$\sigma_i$).
\end{proof}
Now, as we've observed above, the cyclic structures on
$\Omega^R(\hc)$, induced by Karoubi operator and that on
$\tilde\hc_{(\epsilon,1)}^\sharp$, induced from
$\hc_{(\epsilon,1)}^\sharp$ coincide. In effect, we already know,
that $$\Omega^R_n(\hc)\cong\ker\epsilon^{\otimes
n}=(\tilde\hc_{(\epsilon,1)}^\sharp)_n.$$ The only problem is, that
under this isomorphism, differential $\tilde b'$ on
$\hc^\sharp_{(\epsilon,1)}$ corresponds to differential $d$ on
$\Omega^R(\hc)$ (not to $b'$, or $B'$). Other differentials also
play different roles in these cyclic modules. Indeed, it is easy
to see, that $\tilde B'$ of $\tilde\hc_{(\epsilon,1)}^\sharp$
corresponds to the operator $\sum (\kappa')^jb'$ in $\Omega^R(\hc)$,
while $B'=(\sum\kappa')^jd$ in $\Omega^R(\hc)$.

To cure this problem, recall, (\cite{CQ2}, \S3) that the
super-complex $(\Omega(\hc), b+B)$ is quasi-isomorphic to the
subcomplex $P\Omega(\hc)$, on which
\begin{align}
(\kappa-&1)^2=0,\\ \label{eqr}\kappa = 1&- \frac{1}{n(n+1)}bB,\\
\intertext{and hence} B&=(n+1)d
\end{align}
on $P\Omega_n(\hc)$. Here $P$ is the corresponding projection.
Clearly, the same is true about $\kappa',\ b'$, see formulae
\eqref{eq11"}, \eqref{eq11"'}.

Since the quasi-isomorphism $P$ and chain homotopy $Gd$ (this pair
is called special deformation retraction in \cite{CQ2}, \S3) is
expressed in terms of $\kappa'$ and $d$, we conclude, that
$\Omega^R(\hc)$ is quasi-isomorphic to $$P\Omega^R(\hc)\eqdef
P\Omega(\hc)\bigcap\Omega^R(\hc).$$ So, we see, that homology of
$(\Omega^R(\hc), b', B')$ is equal to the homology of
$(P\Omega^R(\hc), b', {\rm deg}\circ d)$. Here ${\rm deg}$ is the
operator, which multiplies the degree $n$ homogeneous elements by
$n$.

On the other hand, consider the map $Gb'$. It is easy to see, that
pair $(P, Gb')$ verifies all the properties of special deformation
retraction for the $(\Omega^R(\hc), d, \sum(\kappa')^jb')$. Recall,
that the operator $\sum(\kappa')^jb'$ is the image of $\tilde B'$
under the above isomorphism. In fact one just repeats the
reasoning from \cite{CQ2}, p.391. So, we conclude this time, that
periodic cohomology of $(\Omega^R(\hc), d, \sum(\kappa')^jb')$ equals
the periodic cohomology of $P\Omega^R(\hc)$ with induced
differentials. And from \eqref{eqr} it follows, that
$\sum(\kappa')^jb'={\rm deg}\circ b'$ on  $P\Omega^R(\hc)$.

Now it is enough to observe, that in both cases we obtain the
periodic cohomology of mixed complex $(P\Omega^R(\hc), d, b')$
(which is the cohomology of the super-complex
$(P\Omega(\hc),d+b')$).
\end{proof}
By a slight modification of these reasoning we obtain the
following
\begin{cor}
Periodic cohomology of the mixed complex $\Omega^R(\hc)$ with
differentials, induced from $b$ and $B$, is isomorphic to the
periodic Hopf-type cohomologie of \hc.
\end{cor}
\begin{proof}
Just observe, that super complex $(\Omega^R(\hc), b+B)$ is as
before quasi-isomorphic to $(P\Omega^R(\hc), d+b)$, and from
\eqref{eq11"'} and \eqref{eqr} we conclude, that at this
subcomplex $b'=b$.
\end{proof}

 Now note, that,
since the space of left-coinvariants in $\Omega(\hc)$ is also closed under
the Hochschild boundary $b$ and Karoubi operator $\kappa$, one can consider
the corresponding mixed subcomplex and its periodic cohomologies.
\begin{prop}
The antypode $S$ of the Hopf algebra $\hc$ induces an isomorphism
of periodic complex of the mixed complex $(\Omega^L(\hc), b, B)$ to
the periodic complex of $(\Omega^R(\hc), b', B')$.
\end{prop}
\begin{proof}
Observe, that, in virtue of the universal properties of
$\Omega(\hc)$, antipode $S$ can be extended to an
anti-automorphism of $\Omega(\hc)$. Since $S^2=1$ in \hc, the same
equation holds for this extension. Hence we get an involutive
anti-automorphism of the universal differential calculus of \hc.
Now a straightforward computation shows, that this map intertwines
the right and left \hc\/-comodule structures, and differentials
$b, \kappa$ and $b', \kappa'$ in the mixed complex.
\end{proof}
\begin{cor}
Periodic cohomology of $(\Omega^L(\hc), b, B)$ is canonically isomorphic to the periodic
Hopf-type cohomology of \hc.
\end{cor}
\begin{rem}
Note, that the universality property of $\Omega(\hc)$ implies
that, in fact, all the maps, defined at the level of \hc\ can be
extended to this differential calculus. Thus, one can introduce
the structure of differential graded Hopf algebra on
$\Omega(\hc)$.
\end{rem}

\section{General case: arbitrary $\delta$ and $\sigma$}
In this section we shall investigate the case of a general modular
pair in involution $(\delta, \sigma)$. We shall reduce this case
to a variant of the construction, we've just considered.

First recall, that for any character $\xi $ of a Hopf algebra
\hc\ one can introduce the following endomorphism (and even
automorphism) of \hc:
\begin{equation}
\tilde\xi:\hc\to\hc,\quad \tilde\xi(a)=a\star\xi
\eqdef\sum a_{(1)}\xi(a_{(2)}).
\end{equation}
The inverse of $\tilde\xi$ is given by the right convolution with
$\xi^{-1}\eqdef\xi\circ S$.

The map $\tilde\xi$ is, in an evident way, a morphism of algebras,
but it does not respect the coalgebra structure on \hc. That is
$\Delta(\tilde\xi(a))\neq\sum
\tilde\xi(a_{(1)})\otimes\tilde\xi(a_{(2)}).$ In fact, this
equation is substituted for the following two:
\begin{align}
\label{eq3.0} \Delta(\tilde\xi(a))=&\sum
a_{(1)}\otimes\tilde\xi(a_{(2)}),\\ \intertext{and}
\label{eq3.1}\Delta(\tilde\xi(a))=&\sum\tilde\xi(a_{(1)})\otimes
Ad_{\xi}(a_{(2)}),\\ \intertext{where} Ad_{\xi}(a)=&\xi^{-1}\star
a\star\xi=\sum\xi(S(a_{(1)}))a_{(2)}\xi(a_{(3)}).
\end{align}
It is easy to see, that $Ad$ defines an action of the group of
characters of \hc\ on \hc\ by Hopf algebra homomorphisms. One
calls it the {\em  adjoint action\/}.

Since $\Omega(\hc)$ is universal differential calculus, we
conclude, that homomorphism $\tilde\xi$ can be extended to higher
degree forms. By abuse of notation we shall denote this map by the
same symbol $\tilde\xi$. Remark, that equations \eqref{eq3.0} and
\eqref{eq3.1} are fulfilled, in a slightly different form, for
this new map, too. Namely:
\begin{equation}
\label{eq3.2}
\Delta_R(\tilde\xi(\omega))=(Id\otimes\tilde\xi)\Delta_R(\omega)=(\tilde\xi\otimes
Ad_\xi)\Delta_R(\omega),\quad\omega\in\Omega(\hc).
\end{equation}
Meanwhile the left coaction remains unchanged:
\begin{equation}
\label{eq3.3}
\Delta_L(\tilde\xi(\omega))=(Id\otimes\tilde\xi)\Delta_L(\omega).
\end{equation}
In fact, the formulae \eqref{eq3.2} and \eqref{eq3.3} are
particular cases of the following observation. As it is remarked
above, $\Omega(\hc)$ is a differential graded Hopf algebra. Its
diagonal map we shall denote by $\tilde\Delta$ If once again by
abuse of notation, $Ad_\xi$ denotes the automorphism of
$\Omega(\hc)$, induced by the appropriate automorphism of \hc,
then the formulas \eqref{eq3.0} and \eqref{eq3.1} hold with
$\Delta$ substituted for $\tilde\Delta$.

Note, that, since $\tilde\xi$ is a map of differential graded
algebras, one can use it to define a new differential structure on
$\Omega(\hc)$. Namely, put
\begin{equation}
d_\xi(\omega)\eqdef d(\tilde\xi(\omega))=\tilde\xi(d\omega)
\end{equation}
One easily checks the following statement, compare \cite{CQ2}:
\begin{prop}
\begin{description}
\item[{\rm({\em i\/})}]{Differential $d_\xi$ verifies the following equation
\begin{equation}
\label{eq3.4}
d_\xi(\omega_1\omega_2)=d_\xi(\omega_1)\tilde\xi(\omega_2)+(-1)^{|\omega_1|}\tilde\xi(\omega_1)d_\xi(\omega_2).
\end{equation}
}
\item[{\rm({\em ii\/})}]{
Algebra $\Omega(\hc)$, equipped with the differential $d_\xi$ is
the universal example of $\xi-$differential calculi on \hc, that
is, of such graded algebras $\Omega$,that
\begin{enumerate}
\item{$\Omega_0=\hc$;}
\item{$\Omega$ is equipped with a degree $1$ map $d_\Omega,\ d_\Omega^2=0$, called differential;}
\item{automorphism $\tilde\xi$ of \hc\ extends to a degree $0$ automorphism of $\Omega$, commuting with $d_\Omega$;}
\item{its differential $d_\Omega$ verifies \eqref{eq3.4}.}
\end{enumerate}
}
\item[{\rm({\em iii\/})}]{
Any element $\theta$ in $\Omega_n(\hc)$ can in a unique way be
represented in the form
\begin{equation}
\label{eq3.5} \theta=\sum_i a_0^i d_\xi(a_1^i)d_\xi(a_2^i)\dots
d_\xi(a_n^i),
\end{equation}
for some $a_i^\alpha\in\hc$.
}
\end{description}
\end{prop}
All this is checked by a straightforward inspection of
definitions. Below we will denote the universal differential
calculus $\Omega(\hc)$ with differential $d_\xi$ by
$\Omega_\xi(\hc)$. We shall also use the presentation of part
({\em iii\/}) to write down the elements of $\Omega_\xi(\hc)$.

Now it is natural to write down the left and right coactions of
\hc\ on $\Omega(\hc)$ in terms of the formula \eqref{eq3.5}. By
virtue of the formulae \eqref{eq3.2} and \eqref{eq3.3}, one gets
\begin{align}
\Delta_R(d_\xi(\omega))=&(d_\xi\otimes Ad_\xi)\Delta_R(\omega),\\
\Delta_L(d_\xi(\omega))=&(Id\otimes d_\xi)\Delta_L(\omega).
\end{align}
Hence, formulae \eqref{eq5} and \eqref{eq6} become
\begin{gather}
\label{eq3.6} \Delta_R(\theta)=\sum_i a_{0,(1)}^i d_\xi
a_{1,(1)}^i d_\xi a_{2,(1)}^i\dots d_\xi a_{n,(1)}^i\otimes
Ad_\xi(a_{0,(2)}^i a_{1,(2)}^i a_{2,(2)}^i\dots a_{n,(2)}^i),\\
\label{eq3.7}\Delta_L(\theta)=\sum_i a_{0,(1)}^i a_{1,(1)}^i
a_{2,(1)}^i\dots a_{n,(1)}^i\otimes a_{0,(2)}^i d_\xi a_{1,(2)}^i
d_\xi a_{2,(2)}^i\dots d_\xi a_{n,(2)}^i.
\end{gather}
Here we've used the fact, that $Ad_\xi$ is a Hopf algebra
homomorphism.

To put short the above considerations, one can say, that one can
consider the universal $\xi-$differential algebra
$(\Omega_\xi(\hc), d_\xi)$, which consists of linear combinations
of elements of the form $a_0d_\xi a_1 d_\xi a_2\dots d_\xi a_n$,
and on which the Hopf algebra \hc\ coacts on both sides by
formulae \eqref{eq3.6} and \eqref{eq3.7}. We shall use this
notation below, though it is not absolutely necessary, since it is
just another way to speak about the universal calculus
$\Omega(\hc)$.

As before, one can consider the spaces of right- and
left-coinvariants in $\Omega_\xi(\hc)$. For instance, the space of
right ones, $\Omega^R_\xi(\hc)$ consists of the tensor powers of
the space spanned by elements
\begin{multline}
\pi^R(a)=
d(a_{(1)})S(a_{(2)})=d_\xi(a_{(1)})\tilde\xi^{-1}(a_{(2)})S(a_{(3)})=\\
d_\xi(a_{(1)})S_{\xi^{-1}}(a_{(2)}),\quad a\in\ker\epsilon.
\end{multline}
As before, $\Omega^R_\xi(\hc)$ is a d.g. subalgebra in
$\Omega_\xi(\hc)$. This follows directly from \eqref{eq3.6}.

In addition to the usual coinvariants, one can consider the space
of elements $\theta$, such that
\begin{equation}
\Delta_R(\theta)=\theta\otimes\sigma
\end{equation}
for some group-like element $\sigma$. We shall call such elements
(right) $\sigma-$coinvariants. Let $\Omega^R_\sigma(\hc)$
(respectively $\Omega^R_{\xi, \sigma}(\hc)$) denote the space of
right $\sigma-$coinvariants in $\Omega(\hc)$ (resp. in
$\Omega_\xi(\hc)$).

Clearly, since $\Delta_R$ commutes with differential $d$,
$\Omega^R_\sigma(\hc)$ is d.g. subalgebra in $\Omega(\hc)$.
Similar statement holds for $d_\xi$, $\Omega^R_{\xi, \sigma}(\hc)$
and $\Omega_\xi(\hc)$.
\begin{prop}
\label{prop1} The differential $d_\xi$ maps the space of
$\sigma-$coinvariants into itself. Moreover, the space
$\Omega^R_{\xi, \sigma}(\hc)$ is a differential graded
$\Omega^R_\xi(\hc)-$ sub-bimodule in $\Omega_\xi(\hc).$
\end{prop}
\begin{proof}
Note, that the right multiplication by $\sigma$ establishes an
isomorphism between the space of (right) coinvariants and the
space of (right) $\sigma-$coinvariants. The inverse is given by
the multiplication by $\sigma^{-1}\eqdef S(\sigma)$. Hence, any
element in $\Omega^R_{\xi, \sigma}(\hc)$ is representible in the
form
\begin{equation}\label{eq3.8}\theta=\theta'\cdot\sigma,
\end{equation}
for a suitable $\theta'\in \Omega^R_\xi(\hc)$. Hence, it is enough
to show, that $d_\xi(\sigma)\in\Omega^R_{\xi,\sigma}(\hc)$. We
compute:
\begin{equation}
d_\xi(\sigma)=d(\tilde\xi(\sigma))=d(\sigma\xi(\sigma)),
\end{equation}
Which is, clearly, (right) $\sigma-$coinvariant, since $d$
commutes with coaction. Here we've used the fact, that $\sigma$ is
group-like, i.e. $\Delta(\sigma)=\sigma\otimes\sigma$.

Finally, the fact that $\Omega^R_{\xi,
\sigma}(\hc)=\Omega^R_\sigma(\hc)$ is a left
$\Omega^R_\xi(\hc)-$module is a consequence of the presentation
\eqref{eq3.8}. Since the left multiplication by $\sigma^{\pm1}$
also establishes an isomorphism between $\Omega^R_{\xi,
\sigma}(\hc)$ and $\Omega^R_\xi(\hc)$, the conclusion follows.
\end{proof}

Let's now define the $\xi-$twisted cyclic structure on
$\Omega_\xi(\hc)=\Omega(\hc)$, that is the analogs of Hochschild
operator $b$ (or $b'$) and Karoubi operator $\kappa$ (or
$\kappa'$). In other words, let's use the presentation
\eqref{eq3.5} to define the following operators on $\Omega(\hc)$.
Put (compare \eqref{eqnew1}-\eqref{eq10'})
\begin{align}
b_\xi(\omega\,d_\xi a)&=(-1)^{|\omega|}(\omega\tilde\xi(a) - a\omega),\\
\label{eq3.9}
\kappa_\xi&=1-b_\xi d - d b_\xi,\\
\intertext{or explicitly}
\label{eq3.10}
\kappa_\xi(\omega d_\xi a)&=(-1)^{|\omega|}da\,\omega.\\
\intertext{And, similarly}
b'_\xi(d_\xi a\,\omega)&=\tilde\xi(a)\omega - \omega a,\\
\kappa'_\xi(d_\xi a\,\omega)&=(-1)^{|\omega|}\omega da.
\end{align}
It is clear, that operators $b_\xi$ and $b'_\xi$ are well defind, since $\tilde\xi(1)=1$.
Also observe, that
\begin{align}
\kappa'_\xi&=\kappa_\xi^{-1},\\
b'_\xi&=b_\xi\kappa'_\xi
\end{align}
Once again,
one easily checks that these operators verify all the properties of the
standard ones, listed in \cite{CQ2}, \S3, only few modifications should be
made. In fact, the following
proposition holds (compare \cite{CQ2}, \S3).
\begin{prop}$\quad$
\begin{description}
\item[{\rm ({\em i\/})}]{
\begin{equation*}
b_\xi^2 = (b'_\xi)^2=0.
\end{equation*}
}
\item[{\rm ({\em ii\/})}]{ Following operators commute
\begin{align*}
[b_\xi,\kappa_\xi]=[d,\kappa_\xi]=[d_\xi,\kappa_\xi]&=0,\\
[\tilde\xi,\kappa_\xi]=[\tilde\xi,b_\xi]&=0
\end{align*}
}
\end{description}
Moreover, on elements of $(\Omega_\xi)_n(\hc)$ one has the following
identities:
\begin{description}
\item[{\rm ({\em iii\/})}]$\kappa_\xi^{n+1}d_\xi=\tilde\xi^{-1}\,d_\xi=d$
\item[{\rm ({\em iv\/})}]$\kappa_\xi^n=\tilde\xi^{-1}+ b_\xi\kappa_\xi^n d$.
\item[{\rm ({\em v\/})}]$\kappa_\xi^n b_\xi=\tilde\xi^{-1}\,b_\xi$.
\item[{\rm ({\em vi\/})}]$\kappa_\xi^{n+1}=\tilde\xi^{-1}(1-d b_\xi)$.
\item[{\rm ({\em
vii\/})}]$(\kappa_\xi^n-\tilde\xi^{-1})(\kappa_\xi^{n+1}-\tilde\xi^{-1})=0$.
\item[{\rm ({\em viii\/})}]{Let
\begin{equation}
B_\xi=\sum_{j=0}^n\kappa_\xi^j d_\xi,
\end{equation}
then $B_\xi d_\xi=d_\xi B_\xi =B^2_\xi=0$.}
\item[{\rm ({\em ix\/})}]$\kappa_\xi^{n(n+1)}-1=b_\xi B_\xi=-B_\xi b_\xi$.
\end{description}
\end{prop}
\begin{proof}
Part ({\em i\/}) is checked by a direct inspection of formulas. Part ({\em
ii\/}) follows from part ({\em i\/}), \eqref{eq3.9} and the fact, that
$\tilde\xi$ is a d.g. algebra automorphism, and hence it commutes with
$b_\xi$ and $d_\xi$ (and consequently with $\kappa_\xi$, too). All the rest
is obtained by mimicking the reasoning of the cited paper, taking in consideration the fact, that
$\tilde\xi$ commutes with all the operators, introduced above. For instance:
let's prove part ({\em iv\/}). We compute, using formula \eqref{eq3.10} and the definitions of $b_\xi$ and $d_\xi$:
\begin{align*}
\kappa_\xi^n(a_0d_\xi a_1\dots d_\xi a_n)&=da_1\dots da_n\,a_0\\
&=\tilde\xi^{-1}(d_\xi a_1\dots d_\xi a_n\,\tilde\xi(a_0))\\
&=\tilde\xi^{-1}(a_0 d_\xi a_1\dots d_\xi a_n +(-1)^{n} b_\xi(d_\xi a_1\dots d_\xi a_n d_\xi a_0))\\
&=\tilde\xi^{-1}\bigl(a_0 d_\xi a_1\dots d_\xi a_n +(-1)^n b_\xi\kappa_\xi^n(d_\xi a_0 d_\xi\tilde\xi(a_1)\dots
d_\xi\tilde\xi(a_n))\bigr)\\
&=\tilde\xi^{-1}\bigl(a_0 d_\xi a_1\dots d_\xi a_n + b_\xi\kappa_\xi^nd_\xi (a_0 d_\xi a_1\dots
d_\xi a_n)\bigr)\\
&=(\tilde\xi^{-1}+b_\xi\kappa_\xi^n d)(a_0d_\xi a_1\dots d_\xi a_n),
\end{align*}
since $d_\xi=d\,\tilde\xi$.
\end{proof}

Now we come to the main result of this paper. The following theorem is a straightforward generalization of the
Theorem \ref{theo1} of the section $2$.
\begin{theorem}
Let $(\delta, \sigma)$ be a modular pair in involution. Let
$\xi=\delta^{-1}$. Then the space $\Omega^R_{\xi, \sigma}(\hc)$ of
$\sigma-$coinvariants in $\Omega_\xi(\hc)$ is stable under the Hochschild
and Karoubi operators and periodic cohomology of the induced mixed complex
is naturally isomorphic to the Hopf-type periodic cohomology $HP_{\delta,\sigma}(\hc)$
of A.~Connes and H.~Moscovici.
\end{theorem}
\begin{proof}
is obtained in a way, absolutely similar to the proof of Theorem
\ref{theo1}. First of all, we establish the first part of this statement (once again we prefer to work with primed versions
of cyclic operatoes).

Namely, let's check, that $b'_\xi(\omega)\in\Omega_\sigma(\hc)$ for all
$\omega\in\Omega_\sigma(\hc)$ (compare \eqref{eq12}). Recall, that
$\xi=\delta^{-1}$:
\begin{equation}
\begin{split}
b'_\xi(\omega)&=b'_\xi\bigl(\pi^R(a_1)\otimes\dots\otimes\pi^R(a_n)\sigma\bigr)\\
&=b'_\xi(d_\xi(a_{1,(1)})S_\delta(a_{1,(2)})\omega'\sigma)=\tilde\xi(a_{1,(1)})S_\delta(a_{1,(2)})\omega'\sigma-S_\delta(a_{1,(2)})\omega'\sigma
a_{1,(1)}\\
&=-\Bigl(\pi^R(S(a_{1,(n+1)})a_2)\otimes\dots\otimes\pi^R(S(a_{1,(3)})a_n)\Bigr)\,(S_\delta(a_{1,(2)})\sigma
a_{1,(1)})\\
&=-\Bigl(\pi^R(S(a_{1,(n+1)})a_2)\otimes\dots\otimes\pi^R(S(a_{1,(3)})a_n)\Bigr)\,(S_\delta(a_{1,(2)})S_\delta^2(a_{1,(1)}))\sigma\\
&=-\Bigl(\pi^R(S(a_{1,(n+1)})a_2)\otimes\dots\otimes\pi^R(S(a_{1,(3)})a_n)\Bigr)\,\Bigl(S_\delta(S_\delta(a_{1,(1)})a_{1,(2)})\Bigr)\sigma\\
&=-\pi^R(S(a_{1,(n-1)})a_2)\otimes\pi^R(S(a_{1,(n)})a_3)\otimes\dots\otimes\pi^R(S_\delta(a_{1,(1)})a_n)\otimes\sigma.
\end{split}
\end{equation}
We've used the fact, that $S_\delta^2(a)=\sigma a\sigma^{-1}$, and the following properties of
$S_\delta$:
\begin{align}
S_\delta(ab)&=S_\delta(b)S_\delta(a);\\
\Delta(S_\delta(a))&=S(h_{(1)})\otimes S_\delta(h_{(2)});\\
S_\delta(h_{(1)})h_{(2)}=\delta(h).
\end{align}
All this is proven by direct computations (see, e.g. \cite{Crai}).

Similarly to the observation, following the equation \eqref{eq12}, one
concludes, that $\kappa'_\xi$ maps $\Omega_\sigma(\hc)$ to itself by a mere
inspection of definitions. But we prefer to give an explicit proof here, too.
We compute (c.f. \eqref{eq13}):
\begin{equation}
\label{eq3.11}
\begin{split}
&\kappa'_\xi(d_\xi
a_{1,(1)}S_\delta(a_{1,(2)})\omega'\sigma)=(-1)^{|\omega'|}S_\delta(a_{1,(2)}\omega'\sigma\,da_{1,(1)}\\
&=(-1)^{|\omega'|}\Bigl(\pi^R(S(a_{1,(n+1)})a_2)\otimes\dots\otimes\pi^R(S(a_{1,(3)})a_n)\Bigr)\,(S_\delta(a_{1,(2)})\sigma
da_{1,(1)}).
\end{split}
\end{equation}
Now, let's consider the last term of this expression separately
(we omit subscript 1 for the sake of brevity):
\begin{equation}
\label{eq3.12}
\begin{split}
S_\delta(a_{(2)})\sigma da_{(1)}&=S_\delta(a_{(3)})\sigma
\pi^R(a_{(1)})a_{(2)}\\&=\pi^R(S(a_{(4)})\sigma
a_{(1)}-S(a_{(4)})\sigma\epsilon(a_{(1)}))S_\delta(a_{(3)})\sigma
a_{(2)}\\&=\pi^R(S(a_{(4)})\sigma
a_{(1)})S_\delta(a_{(3)})S_\delta^2(a_{(2)})\sigma\\&\qquad\qquad-\pi^R(S(a_{(3)})\sigma)S_\delta(a_{(2)})S_\delta^2(a_{(1)})\sigma\\
&=\pi^R(S(a_{(4)})\sigma
a_{(1)})S_\delta(S_\delta(a_{(2)})a_{(3)})\sigma\\&\qquad\qquad-\pi^R(S(a_{(3)})\sigma)S_\delta(S_\delta(a_{(1)})a_{(2)})\sigma\\
&=\pi^R(\delta(a_{(2)})S(a_{(3)})\sigma a_{(1)})\sigma
-\pi^R(\delta(a_{(1)})S(a_{(2)})\sigma)\sigma\\
&=\pi^R(S_\delta(a_{(2)})
S^2_\delta(a_{(1)})\sigma)\sigma-\pi^R(S_\delta(a)\sigma)\sigma\\
&=-\pi^R((S_\delta(a)-\delta(a))\sigma)\sigma.
\end{split}
\end{equation}
Now, equations \eqref{eq3.11} and \eqref{eq3.12} show, that,
identifying $(\Omega^R_\sigma)_n(\hc)$ with
$(\ker\epsilon)^{\otimes n}$ (see proposition \ref{prop1}), one
can write down the twisted Karoubi operator $\kappa_\xi$ as
follows: $$\kappa_\xi(h_1,h_2,\dots,h_n)={\rm
proj}''S_\delta(h_1)(h_2,\dots,h_n,\sigma),$$ where ${\rm proj}''$
is the following projection
\begin{equation}
{\rm
proj}''(h_1,h_2,\dots,h_n)=(h_1-\epsilon(h_1)\cdot1,h_2-\epsilon(h_2)\cdot1,\dots,h_n-\epsilon(h_n)\sigma).
\end{equation}
The rest of the proof reproduces the reasoning of section 1. To
make the analogy more evident, it is worth noting, that
$\tilde\xi$ acts trivially on $\Omega^R_\sigma(\hc)$, since
$\delta(\sigma)=1$, and $\xi=\delta^{-1}$.
\end{proof}

\section{Conclusions}
Finally, we shall make few remarks, concerning the possible ways
to generalize the Hopf-type cohomology.

First of all, consider the special case, discussed in section 2.
Since both $\Omega^R(\hc)$ and $\Omega^L(\hc)$ are closed under
the mixed complex differentials of $\Omega(\hc)$, we conclude,
that the subspace of bi-invariants is also a sub-mixed complex in
$\Omega(\hc)$. Moreover, this subcomplex is stable under the
involution $S$. The corresponding periodic and dihedral periodic
(co)homologies we shall denote by $HP_{bi, \epsilon, 1}(\hc)$ and
$HD_{\epsilon, 1}(\hc)$ respectively. The same constructions
allows one to define bi-invariant homology in the case of
arbitrary modular pair ${\delta,\sigma}$. If $\sigma=1$, one
can reproduce the dihedral construction, too. What is the analog
of dihedral (co)homology in the case of arbitrary $\sigma$ is not so
evident.

Note, that if the Hopf algebra \hc\ is cocommutative, the spaces
of left- and right-(co)invariants coincide, so we see, that in
this case bi-invariant cohomology is isomorphic to the Hopf-type
one. In a generic case the answer is not clear. Besides this, it
isn't clear, whether it is possible to define this type of
bi-invariant and dihedral homology in a $\Omega(\hc)$ independent
way.

Another important observation is, that in order to define
the twisted cyclic structure on $\Omega(\hc)$ (which is equivalent, up to a
change of basis, to $\Omega_\xi(\hc)$), we didn't really use the fact that
isomorphism $\tilde\xi$ was the convolution with a character of \hc, nor even
did we use the fact, that \hc\ is a Hopf algebra. One can come along the same
very line for any autmorphism $f$ of any algebra \ac, to define
$f-$twisted cyclic operators on its universal differential calculus
$\Omega(\ac)$. One can denote the corresponding cyclic (respectively negative
cyclic, periodic cyclic, etc.) homology by $HC_f(\ac)$ (resp. $HC^-_f(\ac),\
HP_f(\ac)$, etc.). For example, one can take automorphism
$$f:\hc\to\hc,\quad f(a)=\alpha\star a\star\beta$$
($\alpha$ and $\beta$ are characters of \hc). Then, if $(\sigma S_{\alpha,\beta})^2=1$ then, passing to $\sigma-$coinvariants
one obtains the construction of \cite{theil} ($S_{\alpha,\beta}$ is the
evident generalization of the map $S_\delta$).

In fact, the homology $HC_f(\ac)$ can be defined in a quite $\Omega(\ac)$ independent way:
see for example \cite{tuset}. Namely, define the $f-$twisted cyclic module as
follows $CC_n(\ac,\,f)=\ac^{\otimes n+1}$ and \begin{align}
\delta^f_i(h_0,h_1,\dots,h_n)&=\begin{cases}
                     (h_0,\dots,h_ih_{i+1},\dots,h_n),&i=0,dots,n-1,\\
                     (f(h_n)h_0,\dots,h_{n-1}),&i=n,
                             \end{cases}\\
\sigma^f_i(h_0,h_1,\dots,h_n)&=(h_0,\dots,h_{i-1},1,h_i,\dots,h_n),\quad 1\leq
i\leq n,\\
\tau^f_n(h_0,h_1,\dots,h_n)&=(-1)^n(f(h_n),h_0,\dots,h_n).
\end{align}
Then all the usual equations of the cyclic operations are fulfilled for
this ones, save that one should substitute the identity operator for an
appropriate tensor power of $f$ in certaine formulae. Further, one defines
the twisted homology theories in completely usual way, by means of the cyclic
duoble complex.

One more way to generalize the constructions above is to use the
remark in the end of section 2. Namely, traking into consideration
the fact, that $\Omega(\hc)$ is a d.g. Hopf algebra, one can
consider it as the input of Connes-Moscovici construction in the
form presented in this paper. Then the universality property of
$\Omega(\hc)$ guarantees, that $S_\delta$ extends to a
homomorphism of this algebra, in such a way, that all the
properties of this map are valid for the extension, too. What one
obtains in this way, is a construction very similar to the
non-commutative Weil complex of Crainic (\cite{Crai}). On the
other hand a very similar construction was introduced by \Dj ur\dJ
evi\v{c} in the guise of universal characteristic classes
construction of Galois-Hopf extensions. This matters will be a
subject of thorough discussion in a following paper.

Finally, there are two more possible approaches to generalizing
constructions, presented in this paper.

The first one consists of substituting the subcomodue
$\Omega_\sigma^R(\hc)$, determined by the modular pair for an
arbitrary cyclically-stable one. Here one can plug in both the
standard and $\xi-$twisted cyclic structures. For instance, if
$\delta=\epsilon$, such stable subcomodules are in one-one
correspondence with all subcoalgebras $\hc'$ of \hc, for which
$$S(a_{(2)})\hc'a_{(1)}\subseteq\hc'$$ for all $a\in\hc$. This is always the case, if \hc\ is commutative. And if \hc\ is
cocommutative, this is equivalent to saying, that $\hc'$ is stable under the
adjoint action of \hc\ on itself.

The second construction seems to be even more general. It consists
of the following idea: it is a well-known fact (see \cite{Kar}),
that one can obtain a good variant of cyclic-type homology, called
{\em the non-commutative De Rham\/} homology from the universal
differential calculus $\Omega(\ac)$ of an algebra \ac\ by passing
to the quotient space
$$\bar\Omega(\ac)\eqdef\Omega(\ac){\big/}[\Omega(\ac),\,\Omega(\ac)],$$
where $[\Omega(\ac),\,\Omega(\ac)]$ is the subspace of graded
commutators of elements of $\Omega(\ac)$. One easily checks, that
the differential $d$ of $\Omega(\ac)$ descends to a differential
in $\bar\Omega(\ac)$. A less trivial fact is, that the induced
homology of $\bar\Omega(\ac)$ coincide with well-defined a
subspace in the cyclic homology of \ac\ (see the original paper of
Karoubi, \cite{Kar}).

Now, if we pass to the barred complex in the case of a Hopf
algebra \hc, we can no more say, that \hc\ acts on it. In fact,
this is not the case, unless \hc\ is commutative. But it is easy
to see, that the space of commutators $[\hc,\,\hc]$ is a coideal
in \hc, hence one can substitute \hc\ for the coalgebra
$\bar\hc=\hc/[\hc,\,\hc]$. Then $\bar\hc$ coacts on
$\bar\Omega(\hc)$ on the right (and on the left, too) and it is
possible to consider the space of coinvariants of this coaction,
namely, the space of those elements
$\bar\omega\in\bar\Omega(\hc)$, which are sent to
$\bar\omega\otimes\bar1$, where $\bar1$ is the group-like element
in $\bar\hc$ determined by $1\in\hc$.

This construction seem to play an important role in the theory of
characteristic classes of Galois-Hopf extensions, which will be an
object of discussion in the next paper. Here we confine ourselves to the
following remark.

An important application of the Hopf-type cohomology is the theory
of characteristic classes of a Hopf-module algebra. On the other
hand, to any Hopf-module algebra one can associate its smashed
product with \hc, which is an example of Galois-Hopf extension of
an algebra. In the next paper we shall investigate the relation
between the Connes and Moscovici construction of characteristic
classes of a Hopf-module algebra and various constructions of
characteristic classes of Galois-Hopf extensions which exist.


\begin{thebibliography}{99}
\bibitem{CM1} Connes~A., Moscovici~H. Hopf algebras, cyclic
cohomology and the transverse index theorem, Communs Math. Phys.
{\bf 198} (1998), 199-246
\bibitem{CM2}Connes~A., Moscovici~H. Cyclic cohomology and Hopf
algebras, Letters Math. Phys. {\bf 48} (1999), 97-108
\bibitem{CM3} Connes~A., Moscovici~H. Cyclic cohomology and Hopf
symmetry, preprint: math.QA/000215
\bibitem{Crai}Crainic~M. Cyclic cohomology of Hopf algebras and a
noncommutative Chern-Weil theory, preprint: math.QA/9812113
\bibitem{CQ1} J.~Cuntz, D.~Quillen. Algebra extensions and
nonsingularity, J. Amer. Math. Soc. {\bf 8}, n.2, (1995), 251-289
\bibitem{CQ2} J.~Cuntz, D.~Quillen. Cyclic homology and
nonsingularity, J. Amer. Math. Soc. {\bf 8}, n.2 (1995), 373-442
\bibitem{Wor} Woronowicz~S.~L. Differential Calculus on Compact Matrix
Pseudogroups (Quantum Groups). Commun. Math. Phys. {\bf 122},
125-170 (1989)
\bibitem{Take} Kreimer~H.~F., Takeuchi~M. Hopf Algebras and Galois Extension
of an Algebra. Indiana Univ. Mathematics Journal, {\bf 30}, (5)
\bibitem{Durd1} {${}$\Dj ur\dJ evi\v{c}~M. Geometry of Quantum Principal Bundles I. Commun.
Math. Phys. {\bf 175}, 457-520 (1996)}
\bibitem{Durd2} {${}$\Dj ur\dJ evi\v{c}~M. Geometry of Quantum Principal Bundles II. Rev.
Math. Phys. {\bf 9}, (5) 531-603 (1997)}
\bibitem{Durd3} ${}$\Dj ur\dJ evi\v{c}~M. Characteristic Classes of Quantum Principal
Bundles. Preprint, Institute of Mathematics, UNAM, Mexico (1995)
\bibitem{Maj-Brz} Majid~S., Brzezinski~T. Quantum Goroup Gauge Theory on
Quantum Spaces. Commun. Math. Phys. {\bf 157}, 591-638 (1993)
\bibitem{Haj} Hajac~P.~M. Strong Connections on Quantum Principal
Bundles. Commun. Math. Phys. {\bf 182}, 579-617 (1996)
\bibitem{Loday} Loday~J-L. Cyclic Homology. A series of Comprehensive Studies
in Mathematics {\bf 301}, Springer-Verlag (1992)
675-692 (1981)
\bibitem{tuset} J.~Kustermans, G.~J.~Murphy, L.~Tuset. Differential
Calculi over Quantum Groups and Twisted Cyclic Cocicles, preprint:
math.QA/0110199 v2
\bibitem{theil}  R.~Taillefer. Cyclic Homology of Hopf Algebras,\\
preprint: math.QA/0009213 v2
\bibitem{iran1} R.~Akbarpour, M.~Khalkhali. Hopf Algebras
Equivariant Cyclic Homology and Cyclic Homology of Prossed Product
Algebras, preprint: math.KT/0011248 v2
\bibitem{iran2} M.~Khalkhali, B.~Rangipour. A New Cyclic Module
for Hopf Algebra, preprint: math.KT/0010153 v2
\bibitem{Sweed} Sweedler~M.~E. Hopf Algebras. W.~A.~Benjamin, Inc.,
New-YorK, 1969
\bibitem{Kar} Karoubi~M. Homologie ciclique et K-th\'eorie.
Ast\'erisque {\bf 149}, 1987
\end{thebibliography}
\end{document}